\documentclass[12pt]{article} %
\pdfoutput=1

\usepackage{amsmath}
\usepackage{amssymb}
\usepackage{latexsym}
\usepackage{bbm}
\usepackage{graphicx}

\newcommand{\mc}[1]{\mathcal{#1}}
\newcommand{\mb}[1]{\mathbf{#1}}
\newcommand{\mbb}[1]{\mathbb{#1}}
\newcommand{\be}{\begin{equation}}
\newcommand{\ee}{\end{equation}}
\newcommand{\bea}{\begin{eqnarray}}
\newcommand{\eea}{\end{eqnarray}}
\newcommand{\ba}{\begin{array}}
\newcommand{\ea}{\end{array}}
\newcommand{\baa}{\left\{ \begin{array}}
\newcommand{\eaa}{\end{array} \right.}
\newcommand{\dt}{\partial_t}
\newcommand{\dx}{\partial_x}
\newcommand{\dy}{\partial_y}
\newcommand{\dz}{\partial_z}
\newcommand{\dzz}{\partial_{zz}}
\newcommand{\dxx}{\partial_{xx}}
\newcommand{\dyy}{\partial_{yy}}

\newtheorem{rmq}{Remark}

\newcommand{\etal}{{\it et al\/}\ }
\newcommand{\dash}{------}

\begin{document}

\title{Variational assimilation of Lagrangian data in oceanography}

\author{Ma\"elle Nodet\\
maelle.nodet@inria.fr}

\maketitle

\begin{abstract}
We consider the assimilation of Lagrangian data into a primitive equations circulation model of the ocean at basin scale. The Lagrangian data are positions of floats drifting at fixed depth. We aim at reconstructing the four-dimensional space-time circulation of the ocean. This problem is solved using the four-dimensional variational technique and the adjoint method. In this problem the control vector is chosen as being the initial state of the dynamical system. The observed variables, namely the positions of the floats, are expressed as a function of the control vector via a nonlinear observation operator. This method has been implemented and has the ability to reconstruct the main patterns of the oceanic circulation. Moreover it is very robust with respect to increase of time-sampling period of observations. We have run many twin experiments in order to analyze the sensitivity of our method to the number of floats, the time-sampling period and the vertical drift level. We compare also the performances of the Lagrangian method to that of the classical Eulerian one. Finally we study the impact of errors on observations.
\end{abstract}


\section{Introduction} 
The world's oceans play a crucial role in governing the earth's weather and climate. Lack of data has been a serious problem in oceanography for a long time. Since ten years the number of observations has greatly increased, with the availability of satellite altimeter data (ie measurements of the free-surface height of the ocean) from Geosat, Topex/Poseidon, Jason and other satellites. In addition to these remote-sensing data we have in situ data, from scientific ships, surface mooring buoys or Lagrangian drifting buoys. Among these observations, Lagrangian data, ie positions of drifting floats, play a relevant role for many reasons: firstly their horizontal coverage is very wide (the whole Atlantic Ocean, for example), secondly they give information about currents, temperature and salinity in depth which are complementary to surface information given by satellite altimeters. For these reasons many national and international programs are organized to deploy drifting floats in the world's oceans. The largest program of this type is Argo (2\,055 floats on the 13th October 2005, 3\,000 planned), whose floats provide also temperature and salinity profiles.\\
There are different types of drifting buoys. In the framework of ocean basin-scale localized experiments, oceanographers have datasets from acoustic floats. These floats emit acoustic signals which are recorded by moored listening stations, and the floats positions are calculated every six hours by triangulation. Large datasets are available especially in the Atlantic Ocean (SAMBA, ARCANE-Eurofloat, ACCE experiments). On a larger scale Argo floats are deployed in order to provide vertical temperature and salinity profiles. Assimilation of Argo thermohaline data has been successfully investigated by Forget (PhD Thesis 2004). Argo floats provide also Lagrangian information, which are their positions every ten days. Indeed they drift freely at a predetermined parking depth (around 1\,000 meters), every ten days they descent to begin profiles from greater depth (2\,000 meters) then they go back to the surface and they record temperature and salinity profiles during ascent. On the surface they transmit data to satellite and they are located by GPS. Thus many different floats networks and Lagrangian datasets are available.\\
In parallel, modeling of the ocean system has greatly improved in both quality and realism, and there are many Ocean Global Circulation Models (OGCM), like for example the OPA PArallelized Ocean model (see Madec \etal 1998). A crucial issue for oceanographers is then to take the best advantage of different types of information included in models to one hand and in various observations to the other hand. Data Assimilation (DA) covers all theoretical and numerical mathematical methods which allow to blend as optimally as possible all sources of information (see the review by Ghil \etal 1997 and by De Mey 1997). There are two main categories of DA methods: variational methods based on optimal control theory (Lions 1968) and statistical ones based on optimal statistical estimation (Jazwinski 1970). Adjoint method is the prototype of variational methods, introduced in meteorology by Penenko and Obraztsov (1976). Its effective implementation in the framework of atmospheric Data Assimilation, namely four dimensional variational assimilation (4D-Var), has been studied by Le Dimet and Talagrand (1986, see also Talagrand and Courtier 1987). Introduction of 4D-Var in oceanography is even more recent (see Thacker and Long 1988, Sheinbaum and Anderson 1990). The prototype of sequential methods is the Kalman filter, introduced in oceanography by Ghil (1989) (see also the reviews by Ghil and Malanotte-Rizzoli 1991).\\
Assimilation of Lagrangian data is in the pipeline. Kamachi and O'Brien (1995) used the adjoint method in a Shallow-Water model with upper-layer thickness as control vector. More recently Mead (2004) has implemented a variational method based on the use of Lagrangian coordinates for Shallow-Water equations. Molcard \etal (2003) and \"Ozg\"okmen \etal (2003) implemented optimal interpolation (which is a simple sequential method) in a reduced-gravity quasi-geostrophic model and in a primitive equations model; their method is based on conversion of Lagrangian data into velocity information. Ide, Kuznetsov, Jones and Salman (2002, 2003, 2005) used Extended and Ensemble Kalman methods to assimilate Lagrangian data into a Shallow-Water model; their method is based on an augmented state vector approach which does not require the conversion of the positions into velocity data. These teams have used simulated data in the twin experiments approach: they don't use real Lagrangian data but idealized observations simulated from a known ``true state'' of the ocean.  Beside these studies Assenbaum and Reverdin (2005) assimilate real data available during the POMME experiment, including Argo floats data, into a very high resolution model thanks to optimal interpolation.\\
Previous works on Lagrangian DA were either based on sequential methods or on variational ones into very simple models. In this paper we investigate variational assimilation of drifters positions into the high resolution primitive equations model OPA.\\
 The aim of variational assimilation methods is to identify the initial state of an evolution problem which minimizes a cost function. This cost function represents the difference between observations and their corresponding model variables. It is minimized using a gradient descent algorithm. The gradient is computed by integration of the adjoint model. Thanks to this formulation there is no need to convert Lagrangian data into velocities data: we can use directly the position observations, although they are not variables of the ocean model, but nonlinear functions of the state variables. Moreover this method is a four dimensional one because the temporal dimension of the observations, ie their Lagrangian nature, is taken into account. The cost function involves a so-called observation operator, which links the state variables (here the velocities) and the observed data (here the positions of drifting particles). This operator is nonlinear and consequently the cost function is not necessarily convex so we used an incremental method (see Courtier \etal 1994) in order to achieve and accelerate the minimization.\\
We implement our method using the primitive equations model OPA. Our configuration is an idealized wind-driven mid-latitude box model, which is representative of the different processes that are going on in the real mid-latitude ocean, as shown by Holland (1978). Then we use the twin experiments approach. As we have said before, Lagrangian observations are simulated from a known ``true state'', so that data are perfectly consistent with the model and it ensures there is no systematic bias in the observations. It is of course unrealistic, but twin experiments are a necessary first-step to validate our method. Indeed in this framework we know exactly the system true state and so we are able to quantify the efficiency of our method by comparing assimilated and true states. Moreover it was relevant not to use real data for many reasons: firstly, Argo floats have not been launched to provide Lagrangian information, the feasibility of exploiting their positions is absolutely not ensured. Secondly, Lagrangian datasets (from Argo to acoustic floats) are very diversified in terms of number of floats, time-sampling period of observations and drifting depth and we want to investigate the sensitivity of our method to these parameters. In order to take into account difficulties of real data (such as drift during ascent and descent for Argo floats or acoustic positioning problems) we also study the impact of errors in observations on assimilation efficiency.\\
The paper is organized as follows: in section \ref{sec:02} we describe the physical model and the Lagrangian simulated data. In section \ref{sec:03} we present the assimilation method and its implementation. Some numerical results are given and commented in section \ref{sec:04}. We conclude in section \ref{sec:5}.

\section{Physical model and Lagrangian data} \label{sec:02}

\subsection{The Primitive Equations of the ocean} 
The ocean circulation model used in our study is a Primitive Equations (PE) model. These equations are derived from Navier Stokes equations (mass conservation and momentum conservation, included the Coriolis force) coupled with a state equation for water density and heat equation, under Boussinesq and hydrostatic approximations (for more details see Lions \etal 1992 and Temam and Ziane 2004).\\
These equations are written as
\be 
\baa{ll} 
\dt u - b \Delta u +(U.\nabla_2)u +w\dz u- a v +  \dx p   =  0  &\textrm{in } \Omega \times (0,t_f)\\
\dt v - b \Delta v +(U.\nabla_2)v +w\dz v + a u +  \dy p  =  0 &\\
 \dz p  - g T  = 0 &\\ 
\partial_{t} T -b \Delta T + (U.\nabla_2) T + w \dz T  +f w =  0  & \textrm{in } \Omega \times (0,t_f)\\
w(x,y,z)=-\int_0^z \dx u(x,y,z') + \dy v(x,y,z') \, dz'& \textrm{in } \Omega \times (0,t_f)\\
U(t=0) = U_0, \qquad T(t=0) = T_0& \textrm{in } \Omega
\eaa  
\ee
where\\
\indent - $\Omega = \Omega_2 \times(0,1)$ is the circulation basin, where $\Omega_2$ is a regular bounded open subset of $\mbb{R}^2$, $x$ and $y$ are the horizontal variables and $z \in (0,1)$ is the vertical one, $(0, t_f)$ is the time interval;\\
\indent - $U=(u,v)$ is the horizontal velocity, $w$ is the vertical velocity, $T$ the temperature and $p$ the pressure;\\
\indent - $U_0=(u_0,v_0)$ and $T_0$ the initial conditions;\\
\indent - $(\nabla_2 .)$ is the horizontal divergence operator and $\Delta = \dxx+\dyy+\dzz$ the 3-D Laplace operator;\\
\indent - $a$, $b$, $f$, $g$ are physical constants.\\
The space boundary conditions are
\be 
\baa{l} 
\ba{llll}\dz u=\tau_u,&\dz v=\tau_v,&T=0  &\textrm{ on } \Gamma_t\\
u=0,&v=0,&T=0 &\textrm{ on }  \partial \Omega \setminus \Gamma_t \ea \\
\int_{z=0}^1 \dx u + \dy v \, dz= 0 \quad \textrm{ in } \Omega_2
\eaa 
\ee
where $\tau=(\tau_u,\tau_v)$ is the stationary wind-forcing, $\partial \Omega$ is the boundary of $\Omega$ and $\Gamma_t = \Omega_2 \times \{z = 1\} $  is its top boundary.

\subsection{Model and configuration} 
\begin{figure}
\begin{center}
   \includegraphics[width=\textwidth]{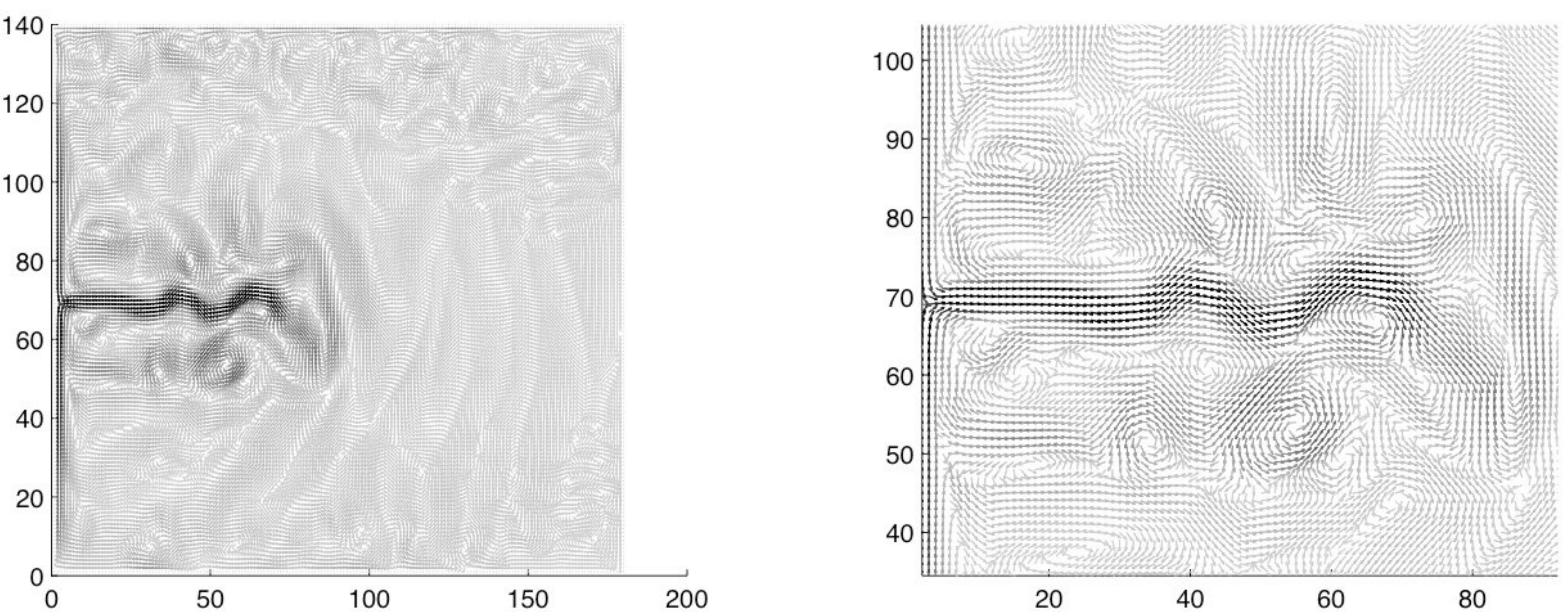}
\end{center}
\caption{\label{fig:01} Instantaneous horizontal velocity field of the true state. On the left the velocity field at the surface on the whole horizontal grid. On the right the velocity field at the surface on a reduced grid centered on the mid-latitude jet. Dark grey vectors represent large velocities.}
\end{figure}

We are using the OPA ocean circulation model developed by LODYC (see Madec \etal 1998), in its 8.1 version.  OPA is a flexible model and can be used either in regional or in global ocean configuration. The prognostic variables are the three-dimensional velocity field $(u,v,w)$ and the thermohaline variables $T$ and $S$. Discretization is based on finite differences in space and time (leap-frog scheme in time). Various physical choices are available to describe ocean physics.\\ 
The characteristics of our configuration are as follows:\\
\indent -- The domain is $\Omega = (0,l)\times (0,L) \times (0,H)$ (longitude, latitude, depth), with $l=2800\textrm{ km}$, $L=3600 \textrm{ km}$ and $H=5000 \textrm{ m}$. It extends from $-56^o$ to $-24^o$ West longitude, and from $22.5^o$ to $47.5^o$ North latitude. \\
\indent -- The horizontal resolution is 20 km, there are 11 vertical levels, so that the number of grid points is 180$\times$140$\times$11 = 277200. \\
\indent -- The time step is $1200$ seconds.\\
\indent -- The model is purely wind-driven.\\
This configuration is a classical eddy-resolving double-gyre circulation. As shown by Holland (1978) this model is representative of real mid-latitude oceans, where circulation is highly nonlinear, non-stationary and where oceanic turbulence is very active. Indeed  a very active and unstable mid-latitude jet develops at the convergence of the subpolar gyre and the subtropical gyre. Non-stationary mesoscale eddies form also along the jet. So this model shows dynamically different processes such as large-scale gyres, mid-latitude jet, mesoscale eddies and also western boundary currents which interact in a complex way. Therefore this configuration is a difficult and interesting situation to study Lagrangian DA.\\
The model is integrated for 25 years until a statistically steady-state is reached, which is our ``true state'' for the twin experiments. Figure \ref{fig:01} shows an instantaneous horizontal velocity field at the surface, on the whole horizontal grid on the left and on a reduced grid on the right. We can see the mid-latitude jet and some mesoscale eddies.

\subsection{Lagrangian data} \label{sec:231}
\begin{figure}
\begin{center}
   \includegraphics[width=\textwidth]{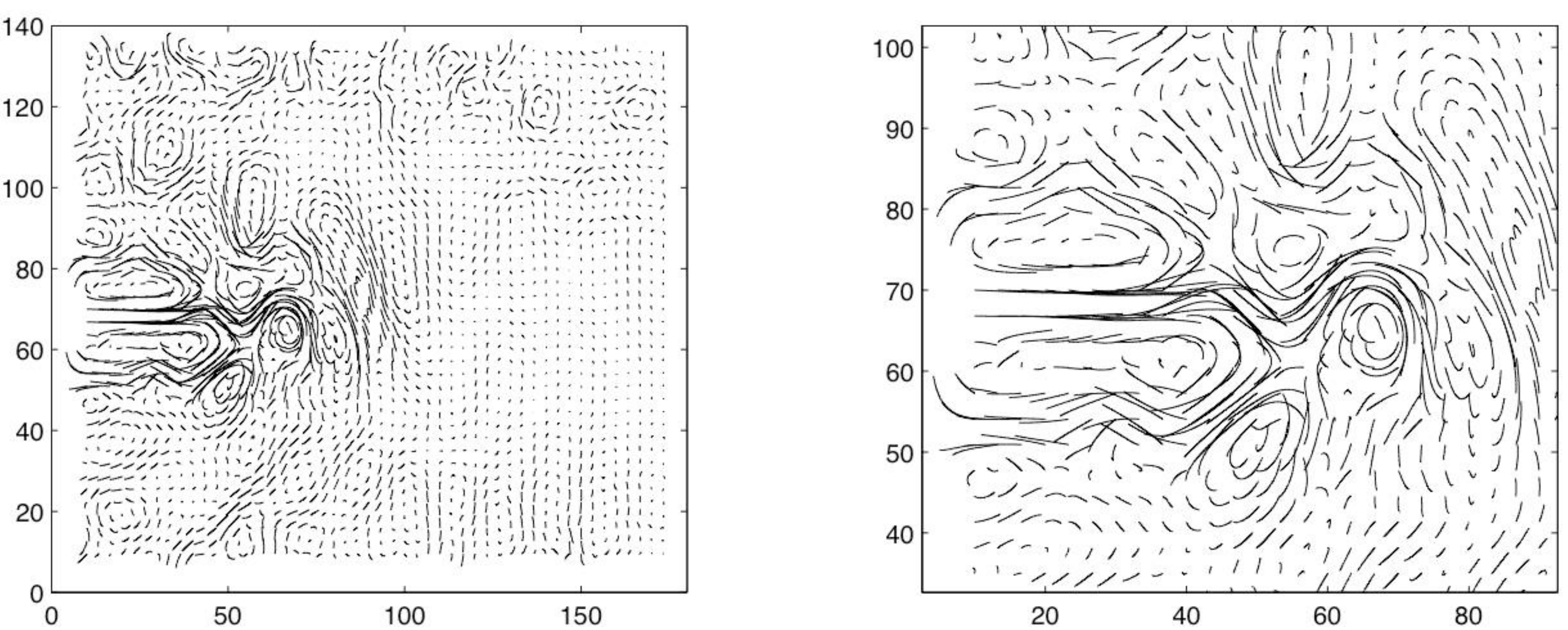}
\end{center}
\caption{\label{fig:02} Trajectories of $2\,000$ floats drifting at level 4 during ten days. On the left, trajectories on the whole horizontal grid at level 4. On the right, trajectories on a reduced grid around the mid-latitude jet. Very different trajectories are observed: short to long, half-circle or straight.}
\end{figure}

Lagrangian data are positions of drifting floats. These floats drift between $z_0-a$ and $z_0+a$ where $z_0$ is given by the user and $a$ is around 25 meters, so that we can assume that the floats drift at fixed depth $z_0$, ie in the horizontal plane $z=z_0$ (Assenbaum 2005, personal communication). 
We denote by $\xi(t) = (\xi_1,\xi_2)(t)$ the position of one float at time  $t$ in the plane $z=z_0$. $\xi(t)$ is the solution of the following differential equation:
\be 
\label{eq:3}
\left\{ \begin{array}{rcl}
\displaystyle \frac{d\xi}{dt}& =& U(t,\xi(t),z_0)\\
 \xi(0)& =& \xi_0 \end{array} \right.
\ee
where $U=(u,v)$ is the horizontal velocity of the flow and $\xi_0$ the initial position of the float. It is important to notice that the mapping $U \mapsto \xi$, which links the variables of the model and the Lagrangian observations is nonlinear.\\
In the twin experiments approach, observations are simulated by the model. The true initial state of the ocean is given and OPA model computes the true velocities of the ocean during a ten-day window. We compute on-line perfect observations.\\
To do that we integrate numerically the equation (\ref{eq:3}) using a leapfrog scheme. This requires the velocity $U$ along the trajectory of the float (ie out of the grid). To achieve this we use the following continuous 2D interpolation '$\textrm{interp}(U,(x,y))$' of the vector field $U$ at the point $(x,y)$:
\[
\ba{l} 
\ba{ll} 
x_1=\lfloor x \rfloor,\qquad &y_1=\lfloor y \rfloor,\\
u_1=U(x_1,y_1),\qquad &u_2=U(x_1+1,y_1),\\
u_3=U(x_1,y_1+1),\qquad &u_4=U(x_1+1,y_1+1),
\ea \\ 
\ba{l}
\textrm{interp}(U,(x,y)) = u_1 + (u_2-u_1)(x-x_1) + (u_3-u_1)(y-y_1) \\
\qquad + (u_1-u_2-u_3+u_4)(x-x_1)(y-y_1)
\ea 
\ea 
\]
where $\lfloor . \rfloor$ denotes the floor function, $(x_1,y_1)$, $(x_1+1,y_1)$, $(x_1,y_1+1)$ and $(x_1+1,y_1+1)$ are the grid points which are the nearest neighbors to $(x,y)$. This function is piecewise affine with respect to $x$ and $y$, continuous with respect to $(x,y)$, linear with respect to $u$. Thus it is not differentiable in $(x,y)$ everywhere. More precisely it is not differentiable at $(x,y)$ if and only if $x=x_1$ or $y=y_1$. It will be a problem to derive the adjoint code, see paragraph \ref{sec:21}. However it is accurate enough to approximate the solution of equation (\ref{eq:3}) and it is very costly to use a differentiable interpolation. Indeed such a method (like cubic splines for example) would compute each interpolated value from the whole field $u$ (ie the values of $u$ at every horizontal grid point) and we would have to inverse a $n$ by $n$ matrix (where $n=25\,200$ is the number of horizontal grid points) at every time step and this is not workable.\\
Let us denote $\xi_k=(\xi_{1,k},\xi_{2,k})$ the horizontal position of the float at time $t_k$, $U$ the horizontal velocity of the fluid at time $t_k$, $U_k$ the velocity at point $\xi_k$, and $h$ the time step of the ocean model. The algorithm step is schematically
\[
\baa{rcl}
\xi_k&=&\xi_{k-2}+2h\,U_{k-1} \\ 
U_{k}&=&\textrm{interp}(U,\xi_k)
\eaa
\]
The dataset is $\{ \xi_{N}, \xi_{2N}, \xi_{3N} ...\}$, where $N$ is an integer. The duration between two data is thus the product of $N$ by the time-step $h$ of the code ($h=1200$ seconds). We call this the time-sampling period.  For example if $N=72$ we have one data per float and per day and the time-sampling period is thus one day.\\
In order to simulate real floats we can add errors to the simulated observations. Origins of errors are multiple: for acoustic floats, inaccuracy can come from acoustic sources (accuracy of their positioning, clocks accuracy, bottom topography -- acoustic shadow problem, etc.), floats (listening period accuracy, complexity of the trajectory, technical problem -- temporary ``deafness'', etc.) or communications quality. For Argo floats errors are due to drift during ascent and descent and also to drift at the surface between ascent/descent and satellite communication. Errors amplitude is around 3 to 4 kilometers for acoustic floats (T. Reynaud, private communication) and 2 to 6 kilometers for Argo floats (M. Assenbaum, private communication).\\
Figure \ref{fig:02} represents perfect data simulated by the algorithm with $2\,000$ floats for 10 days at level 4 (1\,000 meters).

\section{Description and implementation of the variational assimilation method} 
\label{sec:03}

\subsection{Description of the assimilation problem} \label{sec:31}

Without loss of generality we assume that there is only one assimilated data: the dataset is the position $\xi(t_1)=(\xi_1(t_1),\xi_2(t_1))$ in the horizontal plane $z=z_0$ of a single float at a single time $t_1$. We use here the notations established by Ide \etal (1997). We denote by $\mb{y}^o =(\xi_1(t_1),\xi_2(t_1))$ this data.  Our problem is to minimize the following cost function with respect to the control vector $\mb{x}$:
\be 
\ba{lcccc} 
\label{eq:7}
\displaystyle \mc{J}(\mb{x})&=& \frac{1}{2}  \|\mc{G}\mc{M}(\mb{x}) - \mb{y}^o \|^2 &+& \frac{\omega}{2}\|\mb{x}-\mb{x}^b\|_{\mb{B}}^2\\
& =& \mc{J}^o(\mb{x})  &+&\omega \, \mc{J}^b(\mb{x})  
\ea 
\ee
where \\
\indent -- the control vector $\mb{x} = (u_0,v_0,T_0)$ is the initial state vector,\\
\indent -- the $\mb{B}$-norm is calculated thanks to the background error covariance matrix $\mb{B}^{-1}$: $\|\psi\|_{\mb{B}}^2=\psi^T \mb{B}^{-1} \psi$,\\
\indent -- $\mb{x}^b$ is another initial state of the ocean, called the background or first guess, which is required to be close to the minimum $\bar{\mb{x}}$,\\
\indent -- $\mc{M}$ is the discrete ocean model and $\mc{M}(\mb{x})$ is the discrete state vector (one value per variable, per grid-point and per time step),\\
\indent -- $\mc{G}$ is the discrete nonlinear observation operator, which links the state of the fluid (and especially the horizontal velocity $U$) with the data, $\mc{G}\mc{M}(\mb{x})= \xi(t)$ where $\xi(t)$ is defined by equation (\ref{eq:3}),\\
\indent -- $\|.\|$ is the euclidean norm in $\mbb{R}^2$,\\
\indent -- $\mb{y}^o$ is the observations vector.\\
Then $\mc{J}^o$ quantifies the misfit between observations and the state of the system, $\mc{J}^b$ represents the distance (in terms of the $\mb{B}$-norm) between the control vector and the background. It is also a regularization term thanks to which the inverse problem of finding the minimum $\mb{x}^*$ becomes well-posed. The parameter $\omega$ represents the relative weight of the regularization term with respect to the observation term and it must be chosen carefully.

\subsection{Numerical variational assimilation : incremental 4D-Var} 
\label{sec:11}
Four dimensional variational assimilation (4D-Var, see Le Dimet and Talagrand 1986) is an iterative numerical method which aims to approximate the solution  ${\mb{x}^*}$ of discrete assimilation problems with cost function of type (\ref{eq:7}).\\
In 4D-Var a gradient descent algorithm is used to minimize the cost function, the gradient being obtained by solving the discrete adjoint equations. It is an efficient method but it is very costly when the direct model and the observation operator are not linear, for at least two reasons. Firstly every iteration of the adjoint method requires one integration of the full non linear direct model and one integration of the adjoint of the linearized model. Secondly the cost function is not necessarily convex and the minimization process may converge to a local minimum, or it may take considerable time to converge, or may not converge at all.\\
Incremental 4D-Var (see Courtier \etal 1994) avoids, to some extent, both of these problems. In this approach, the nonlinear model is approximated by a simplified linear model (called tangent linear model) and the nonlinear observation operator is linearized around a reference state. The cost function becomes quadratical, it has a unique minimum and this minimum is assumed to be close to the one of the full non quadratical cost function. In that case the minimization process converges quickly. Moreover this approach takes into account weak nonlinearities, because the tangent linear model and the adjoint model are updated three or four times. 
\begin{rmq}\label{rmq:3}
The approximation of the full nonlinear model by the tangent linear model is called the \emph{tangent linear hypothesis} (TLH). In our highly nonlinear configuration, we have to use a ten-day time-window so that the TLH is valid (see section \ref{sec:04}).
\end{rmq}

\subsection{Implementation in OPAVAR} \label{sec:21}

The OPAVAR 8.1 package developed by Weaver \etal (2003) includes the direct non linear model OPA 8.1 developed by LODYC, the tangent linear model, the adjoint model and a minimization module. Weaver has implemented a preconditioning through the $\mb{B}$ matrix, via the change of variables $\delta \mb{w} =\mb{B}^{-1/2}\mb{\delta x}$, following the method introduced by Courtier \etal (1994).
The observation operators of OPAVAR 8.1 are interpolation and projection operators. To assimilate Lagrangian data we have implemented the non linear observation operator (see section \ref{sec:231}), its linearization around the reference trajectory and the adjoint of the linear observation operator.\\
To obtain the tangent (and adjoint) codes of the discrete observation operator we use the recipes for (hand-coding) adjoint code construction of Talagrand (1991) and Giering and Kaminski (1998). The direct and tangent algorithms are schematically:
\begin{itemize}\item Direct code:
\[ 
\baa{l}  \xi_k=\xi_{k-2}+2h\,U_{k-1} \\
 U_k=\textrm{interp}(U,\xi_k)
\eaa 
\]
where 'interp' is the interpolation function of $U$ at point $\xi$ (see section \ref{sec:231}).
\item Linear tangent code: 
\[ 
\baa{l} \delta \xi_k=\delta \xi_{k-2}+2h\,\delta U_{k-1} \\
 \delta U_k=\textrm{interp}(\delta U,\xi_k) + \delta \xi_k . \partial_{(x,y)}\textrm{interp}(U,\xi_k) 
\eaa 
\]
where '$\partial_{(x,y)}\textrm{interp}$' is the derivative of the 'interp' function with respect to $(x,y)$. The term $\partial_{(x,y)}\textrm{interp}$ is specific to Lagrangian data. It leads to a slight difficulty, because the function 'interp' is linear with respect of $U$ but it is not derivable in  $(x,y)$ at points with integer coordinates. Thus we have chosen the values of that derivative at these points, using finite centered differences.
\end{itemize}

\section{Numerical results} \label{sec:04}
In this section we present the results of our numerical experiments. We begin with a brief description of our choices.
\paragraph{Background and time-window width.} In these experiments we have assimilated only Lagrangian data and we assume that the true initial temperature and salinity were known. Background and time-window width are related because of the incremental formulation: indeed the full nonlinear model is linearized \emph{around the background over the whole time-window}. When the background is too different from the true state or the time-window is too wide, approximation errors are large ie the \emph{tangent linear hypothesis} (TLH) is not valid any more. So we compute some correlations to choose both of them. If we denote $\mb{x}^t$ the true state and $\mb{M}$ the tangent linear model, we can compute the nonlinear and linear perturbations $\delta_1$ and $\delta_2$:
\[
\delta_1 = \mc{M}(\mb{x}^t) - \mc{M}(\mb{x}^b),\quad \delta_2 = \mb{M}(\mb{x}^t - \mb{x}^b)
\]
Then we compute (as a function of time) the spatial correlation between $\delta_1$ and $\delta_2$ according to the formula:
\be
\label{eq:1}
Cor(\delta_1,\delta_2) = \frac{\langle \delta_1 \delta_2 \rangle -\langle \delta_1\rangle\langle \delta_2 \rangle }{\sqrt{(\langle \delta_1^2\rangle-\langle \delta_1\rangle^2)(\langle \delta_2^2 \rangle-\langle \delta_2 \rangle^2)}}
\ee
where $\langle X\rangle$ is the spatial mean of $X$. The closer to 1 the correlation is, the better the adequacy between the fields is. Table \ref{tab:6} shows correlation at the end of the time-window between linear and nonlinear perturbations for different time-window widths (10 or 20 days) and different background choices (the state of the ocean 10 days or 1 month before the true initial state). We can see that the TLH is not valid with a 20-day window. In the sequel we use a 10-day time-window and the state of the ocean ten days before the true one as a background state, as in this context the TLH is valid. We ran the model with the background as initial state and we obtained a ``without-assimilation'' state, called background in the sequel. It will be compared to the assimilated state in order to quantify the efficiency of the assimilation process.
\begin{table}
\caption{\label{tab:6}Background and time-window width choices: spatial correlation between nonlinear and linear perturbations at the end of the time-window, according to formula (\ref{eq:1}). The background is the state of the ocean 10 days or 1 month before the true initial state. }
\begin{tabular}{@{}llc@{}}
\hline
time-window width &  background & correlation\\
\hline
10 days & 10 days & 0.80\\
10 days & 1 month & 0.67\\
20 days & 10 days & 0.50\\
20 days & 1 month & 0.42\\
\hline
\end{tabular}
\end{table}

\paragraph{The $\mb{B}$ matrix.} The choice of the $\mb{B}$ matrix is crucial because of its dual purpose (preconditioning and regularization). Firstly we have tested very simple matrices (identity, energy weights) and the results were quite bad: the convergence was very slow and the analysis increments were very noisy. These matrices are used in OPAVAR only for debugging purpose, with extremely idealized assimilation context, for example when the whole state vector is observed ie with data everywhere. In order to obtain smoother increments and to accelerate the convergence, we used then the diffusion filter method (see Weaver and Courtier 2001) which gives good results.

\paragraph{Diagnostics.} Our diagnostics are based on RMS error between the true velocity and the assimilated one, compared with the RMS error between the true velocity and the background one. The RMS error is plotted as a function of time or of the vertical level or of another parameter. For example, we have the following formula for the time-dependent RMS error:
\be 
\label{eq:8}
\textrm{error}(u,t) = \Big( \frac{\sum_{i,j,k} |u_t(i,j,k,t)-u(i,j,k,t)|^2}{\sum_{i,j,k}|u_t(i,j,k,t)|^2} \Big)^{1/2}
\ee
where $u_t$ is the true state, $u$ the assimilated state (or the background), $t$ is the time and $(i,j,k)$ a grid point, where $(i,j)$ are the horizontal coordinates and $k$ the vertical one.\\
We made the following experiments: first experiment and diagnostics in section \ref{sec:41}, sensitivity to the floats network parameters (time sampling of the position measurements, number of floats, drifting level, coupled impact of number and time sampling) in section \ref{sec:42}, comparison with another variational method in section \ref{sec:43}.

\subsection{First experiment}\label{sec:41}
\begin{figure}
\begin{center}
\includegraphics[width=\textwidth]{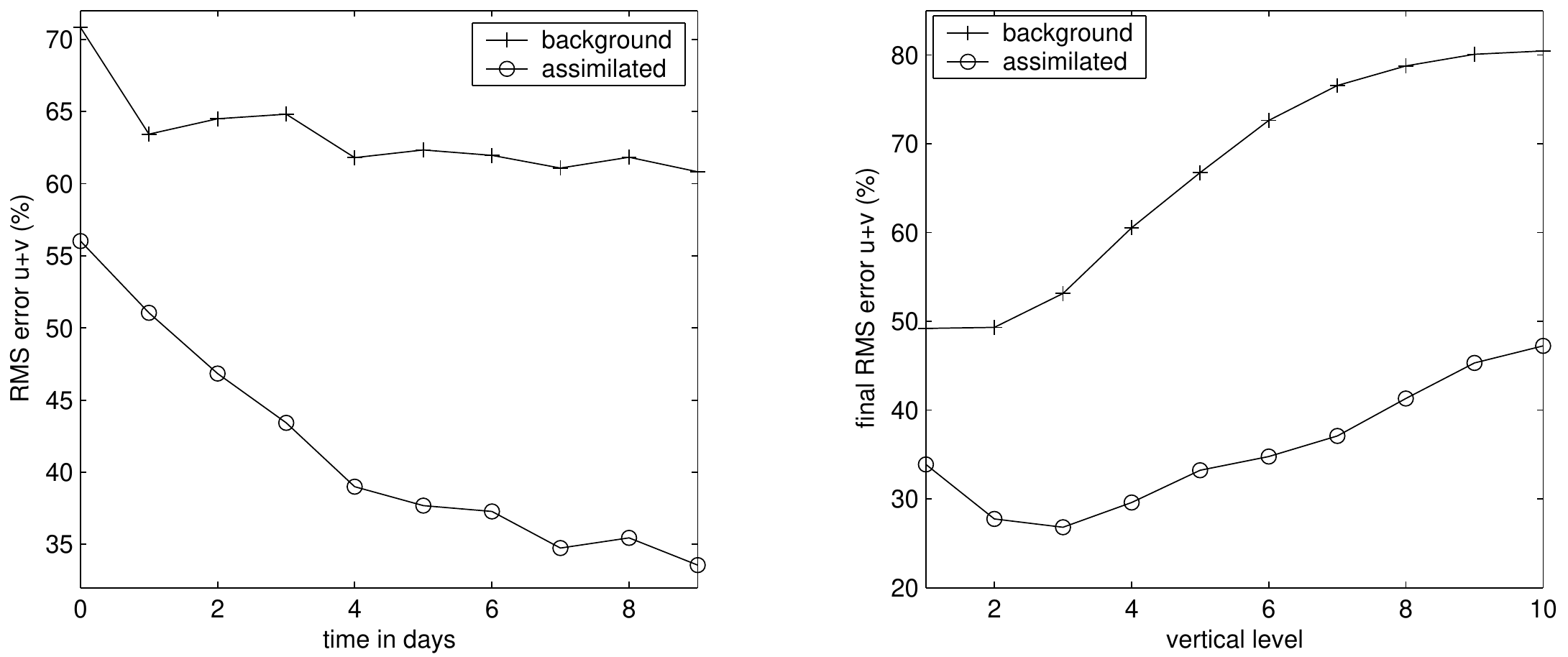}
\end{center}
\caption{ \label{fig:03} First experiment: u+v RMS errors corresponding to the assimilation of the positions (sampled 4 times a day) of 3\,000 floats drifting at level 4. On the left, RMS error as a function of time. On the right, RMS error as a function of the vertical level. For reference the error without assimilation (background) is also displayed.}
\end{figure}

\begin{figure}
\begin{center}
   \includegraphics[width=\textwidth]{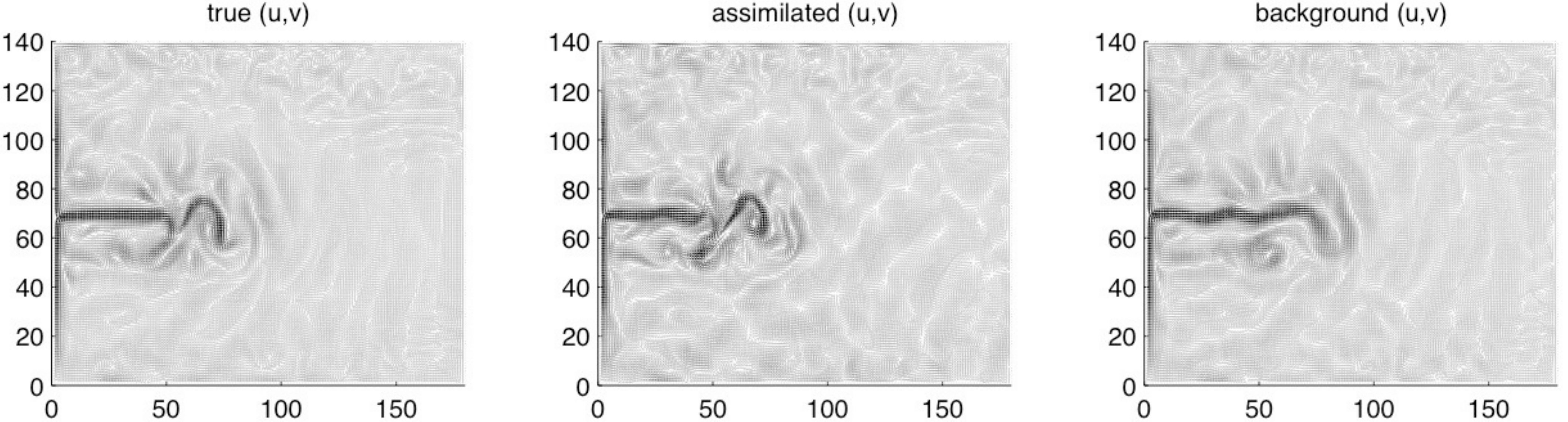}
\end{center}
\caption{ \label{fig:05} First experiment: horizontal surface velocity field at the final time. The true field is displayed on the left. The field corresponding to the assimilation of floats positions is displayed on the middle. For reference the field obtained without assimilation (background) is displayed on the right. }
\end{figure}

We present here the results of a typical experiment. There are $3\,000$ floats drifting at level 4 in the ocean for 10 days. The Lagrangian data  are collected once a day. Thus the total amount of data is $2 \times 3\,000 \times 10 = 60\,000$.\\
Figure \ref{fig:03} (on the left) shows the RMS error of the experiment as function of time, according to formula (\ref{eq:8}). We have put the error for the background (= without assimilation state) on the same plot.\\
Figure \ref{fig:03} (on the right) shows the total RMS error as a function of the vertical level (where 1 represents the surface and 10 the bottom), according to the formula:
\begin{equation}\label{eq:2}
\textrm{error}(u,k) = \Big( \frac{\sum_{i,j,t} |u_t(i,j,k,t)-u(i,j,k,t)|^2}{\sum_{i,j,t}|u_t(i,j,k,t)|^2} \Big)^{1/2}
\end{equation}
We can see that the error with assimilation is twice lower than without. Moreover the assimilation process improves every vertical level and not only the 4th one. 
\begin{rmq} \label{rmq:1}
We can notice that the RMS error at the beginning are quite large. It is explained by the following fact: the relative weight of the regularization term $\mc{J}^b$ (see section \ref{sec:31}) has to be large enough to ensure the convergence of the minimization process, thus the assimilated state is a compromise between background and observations. 
\end{rmq}
Figure \ref{fig:05} shows the horizontal velocity field $U=(u,v)$ at level 1 at the final time. We can notice that the main patterns as the mid-latitude jet and the bigger eddies are quite similar to the true ones.

\paragraph{Longer experiments.} 
We have seen before that the TLH is not valid for a 20-day window, so that we cannot use incremental 4D-Var for longer windows. However we can restart the assimilation process over the next 10-day window in order to do longer experiments: the new background (at the beginning of the new window) is the previous assimilated state (at the end of the previous window), so that the new background carries information from the previous assimilation process. Table \ref{tab:1} shows the relative RMS errors (\ref{eq:2}) for $u+v$ of a 30-day experiment at different times. We can see that error with assimilation is lower than 15\% at the end of the window, ie it is less than one fourth of the error without assimilation. This is a very good result: 3 succesive assimilation processes enable to reconstruct a very good approximation of the true state.
\begin{table}
\caption{\label{tab:1}RMS errors for u+v (in \%) at different instants, during a 30-day experiment, with and without assimilation of the positions of 1\,000 floats drifting at level 4, sampled once a day. }
\begin{tabular}{@{}lllll@{}}
\hline
Experiment & t = 0 & t = 9 days   & t = 19 days & t = 29 days\\
\hline
Without Assimilation &   70.8  &   61.8     &   62.3      &    59.8      \\
\hline
With Assimilation&   55.4   &    34.9    &  21.7      &   13.6      \\
\hline
\end{tabular}
\end{table}

\subsection{Sensitivity to the floats network parameters}\label{sec:42}

In operational oceanography, sensitivity analysis is central to the observational network design. Indeed, in situ and remote observation instruments are very expensive (e.g. one Argo floats costs 15\,000\$) and they must be optimally used. Ngodock (PhD thesis 1996) shows that second order analysis (ie the derivation of the optimality system or \emph{second order adjoint system}, see also Wang \etal 1992) enables to analyze the sensitivity of the 4D-Var assimilation system to the design of the observational network. Second order adjoint information (see also the review paper by Le Dimet \etal 2002) is actually central to adaptive observation network and observation targeting issues. However it requires the storage of model, tangent and adjoint trajectories, so that it is not workable in OPAVAR at the present time because of computer memory limitations.\\

\noindent So we have performed a lot of experiments to analyze the sensitivity to various parameters of our assimilation process. Indeed in the ocean the network's parameters can widely change, from Argo floats (1000-2000 meters depth, one data per 10 days) to acoustic floats (various depth, time sampling period around 6 hours) or drifters in the upper ocean (near surface, time sampling period very short)...\\
Here are the parameters that we consider:\\
\indent - the time sampling period, varying from 6 hours to 10 days,\\
\indent - the number of floats, varying from 300 to 3000,\\
\indent - the vertical level of drift, varying from 1 (surface) to 10 (bottom),\\
We analyze also the coupled effect of the number and the time sampling period.

\subsubsection{Sensitivity to the time sampling period.}
\begin{figure}
\begin{center}
\includegraphics[width=\textwidth]{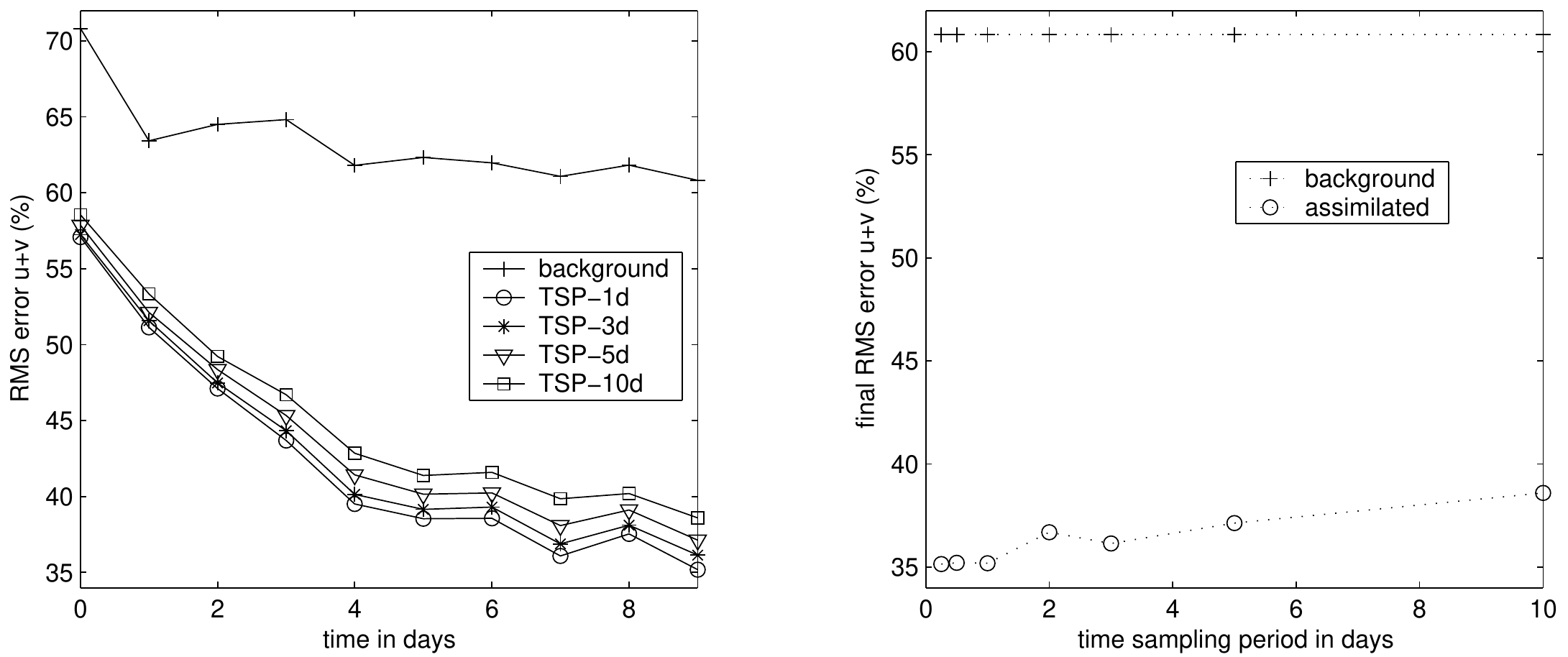}
\end{center}
\caption{\label{fig:06}Sensitivity to the time sampling period: u+v RMS errors corresponding to assimilation of 3\,000 floats' positions with different time-sampling periods of observation. On the left: error as a function of time for each 'TSP' experiment and for the background (without assimilation reference state). On the right: final error as a function of the TSP; the error without assimilation is also displayed.}
\end{figure}

The framework of this experiment is the following: we performed seven different experiments with exactly the same initial conditions, namely $3\,000$ floats at level 4. The only difference in these experiments is the time sampling period, which will be denoted shortly by TSP in the sequel. The experiments are denoted by TSP-xxx where xxx is the time sampling period, in hours (6 or 12h) or in days (1, 2, 3, 5 or 10d). 
Figure \ref{fig:06} shows the RMS error as a function of time for each TSP experiment, except (for readability) experiments TSP-6h, TSP-12h and TSP-2d. It shows also the total RMS error as a function of the time sampling period. We can see that our method is robust with respect to the increase of the time sampling period. This is very encouraging. Indeed every prior study is very sensitive to the TSP and shows quite bad results when the TSP is larger than 2 or 3 days (see Molcard \etal 2003, Mead 2005 and section \ref{sec:43}). Our method does not show this sensitivity, velocities are quite well reconstructed even when the TSP is large and especially with a 10-day period, which is a very positive result. 

\subsubsection{Sensitivity to the number of floats.}
\begin{figure}
\begin{center}
\includegraphics[width=\textwidth]{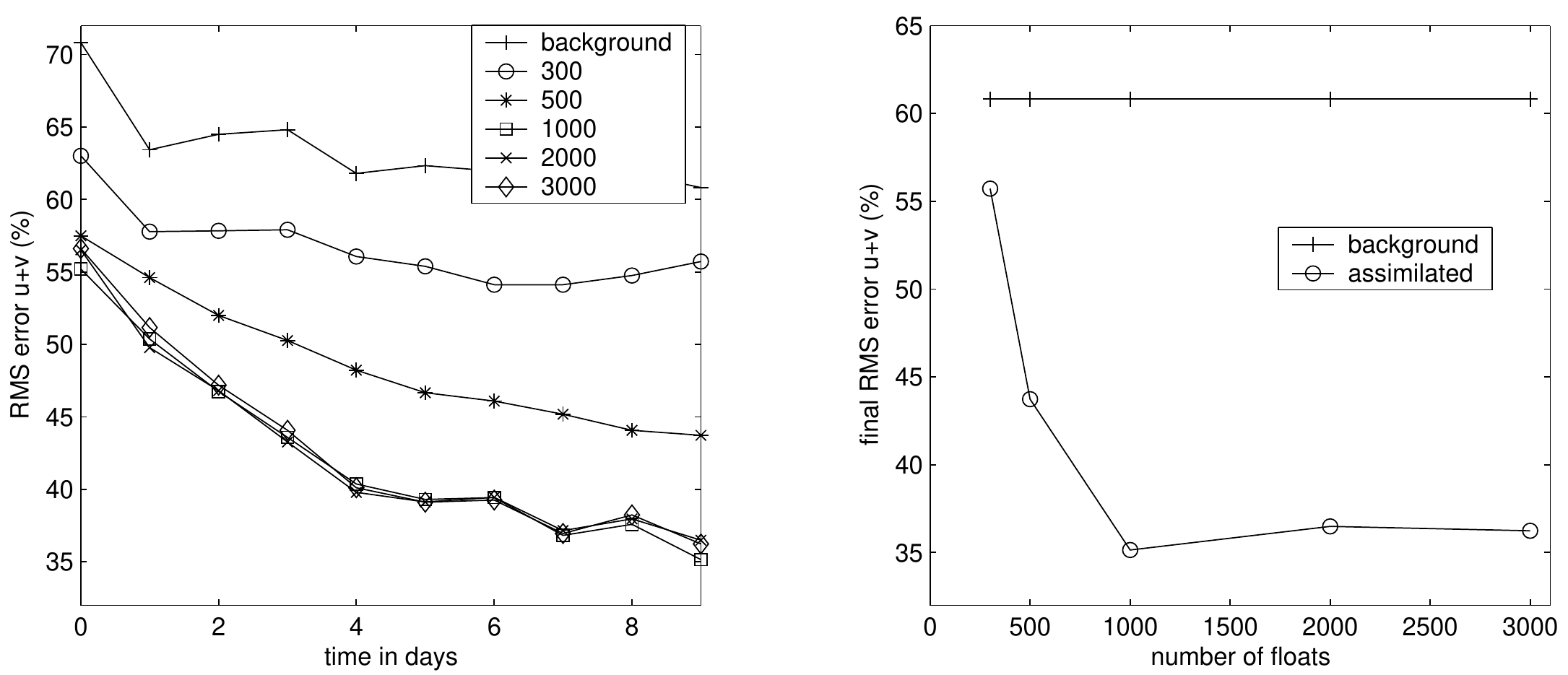}
\end{center}
\caption{\label{fig:07}Sensitivity to the number of floats: u+v RMS errors corresponding to assimilation of 300 to 3\,000 floats' positions. On the left: error as a function of time for each 'NUM' experiment and for the background (without assimilation). On the right: final error as a function of the number of floats; the error without assimilation is also displayed.}
\end{figure}

We perform five experiments with varying number of floats drifting at the same vertical level (4) and with positions sampled with the same period (6 hours). The experiments are denoted by NUM-xxx, where xxx is the number of floats. Figure \ref{fig:07} shows the RMS error as a function of time and the total RMS error as a function of the number of floats for each experiment. We can see that the number of floats has great influence on the results. Under a minimal number (1\,000) the velocities are badly reconstructed, undoubtedly because there is not enough information to constrain the flow. However, when we perform longer experiments, we get satisfactory results for small numbers like 500 and 300. For a ten-day window the results are optimal with a 1\,000 floats network and they don't improve with higher numbers. Obviously the information becomes redundant and it is useless to add floats. The associated density is one float per 10\,000 km$^2$, which is ten times more than the planned Argo density (namely around 100 floats in our configuration). Even if we perform longer experiments, the Argo density is too small to constrain the velocity field. It is more appropriate to use Lagrangian data from localized experiments such as acoustic floats launchings in the Atlantic Ocean (like e.g. SAMBA), whose floats densities are higher.

\subsubsection{Sensitivity to the vertical drift level.}
\begin{figure}
\begin{center}
\includegraphics[width=\textwidth]{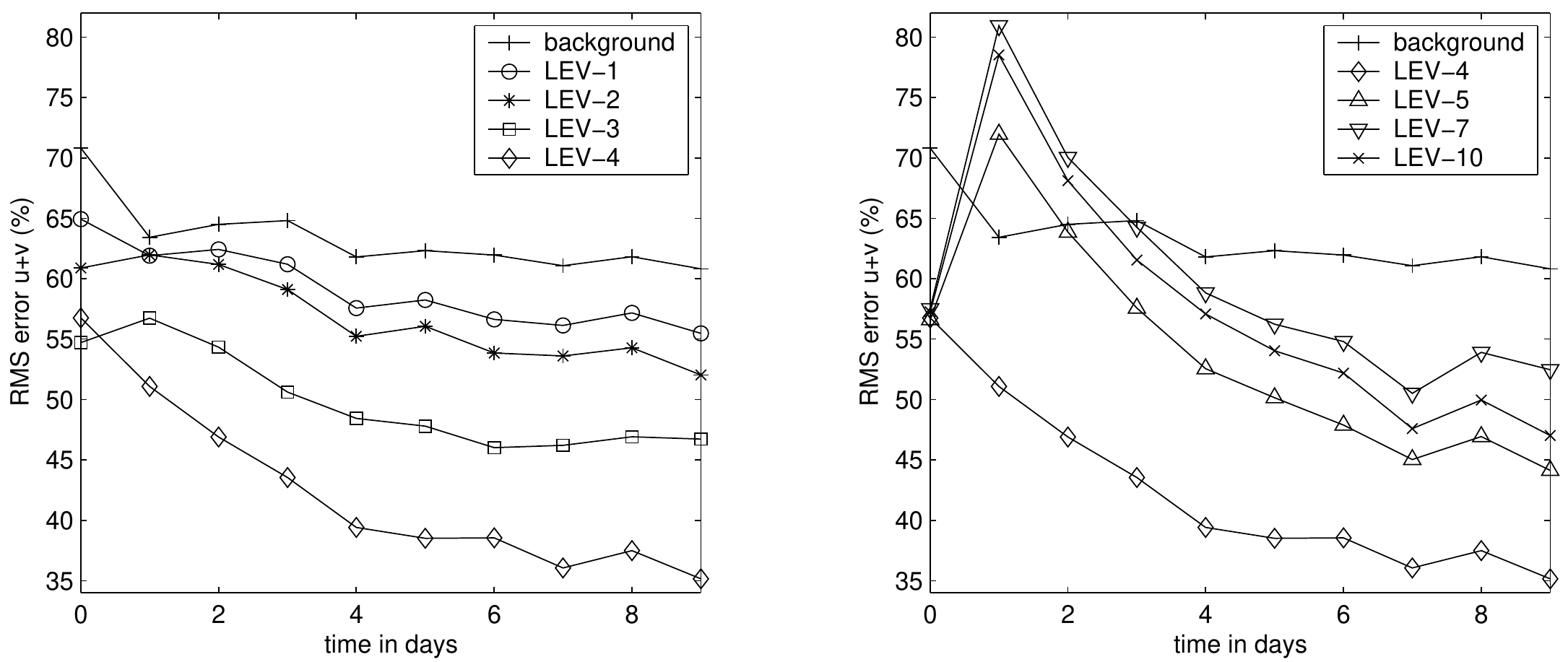}
\end{center}
\caption{\label{fig:08}Sensitivity to the vertical drift level: time-evolution of the u+v RMS errors corresponding to the assimilation of positions of floats drifting at different depths. On the left: error for upper levels experiments (above 1\,000 meters) and for the background (without assimilation). On the right: error for lower levels experiments (below 1\,000 meters) and for the background (without assimilation).}
\end{figure}

Again we perform seven experiments with $3\,000$ floats drifting at varying vertical level and fixed TSP (6 hours). As usually we denote by LEV-x the experiment involving floats at level x. Figure \ref{fig:08} shows the RMS error as a function of time for upper levels on the left and lower levels on the right. Again the results are very sensitive to the position of the floats. The three best levels are 3, 4 and 5, ie the intermediate levels. From a physical point of view it is coherent because the information propagates vertically with a finite velocity so that very upper (1, 2) and very lower (7 to 10) levels are penalized. Moreover upper levels (1 to 4) are the most energetic ones (from the kinetic turbulent energy point of view), quasi ten times more than the lower ones (levels 5 to 10), it seems quite natural that the best results are obtained with floats drifting at level 4 which is both intermediate and energetic.

\subsubsection{Coupled impact of number of floats and time sampling period.}
\begin{figure}
\begin{center}
\includegraphics[width=\textwidth]{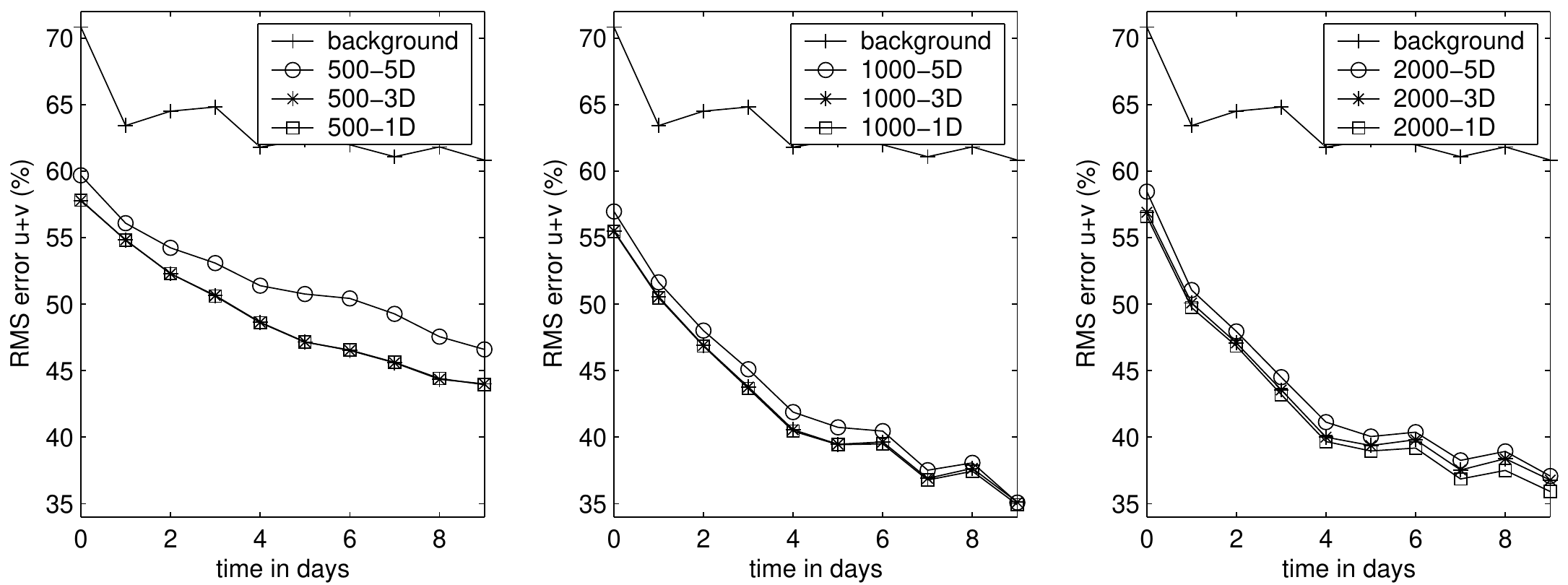}
\end{center}
\caption{\label{fig:09}Coupled impact of number and time sampling period: time-evolution of the u+v RMS errors corresponding to assimilation experiments with 500 to 2\,000 floats and positions sampled every 1 to 5 days: experiment with 500 floats on the left, 1\,000 in the middle and 2\,000 on the right. For reference, the background error is also displayed on each plot.}
\end{figure}

Here we look at the coupled effect of varying number of floats and varying TSP, for example in order to answer the following question: is the total number of data an important variable to measure the efficiency of the assimilation? So we perform nine experiments, denoted by nnn-xxx where nnn is the number of floats and xxx is the TSP. These experiments and their final RMS error are described in Table \ref{tab:4}. Figure \ref{fig:09} represent the RMS error as a function of time for the 500-xxx experiments on the left, 1000-xxx in the middle and 2000-xxx on the right with the same scale on the axis of ordinates. The results are complementary to the precedent experiments. Indeed we can see that 1\,000 is an optimal number for this configuration whatever the TSP and that our method is stable with respect to large TSP whatever the number of floats. Thus we can conclude that, in our configuration, it seems optimal to launch around 1\,000 floats and that the TSP can be chosen quite large.
\begin{table}
\caption{\label{tab:4}Coupled impact of number of floats and TSP: Final RMS error corresponding to assimilation experiments with 500 to 2\,000 floats and positions sampled every 1 to 5 days. Total number of observations is also given. The ``background experiment'' results are also shown for reference.}
\begin{tabular}{@{}ccc@{}}
\hline
Experiment &  Total number of data  & Final Error (\%)\\
\hline
500-5D   &  1\,000  &  46.6  \\
500-3D   &  1\,500  &  44.0  \\
500-1D   &  5\,000  &  44.0  \\
2000-5D  &  4\,000  &  37.1  \\
2000-3D  &  6\,000  &  36.8  \\
2000-1D  &  20\,000 &  35.9  \\
1000-5D  &  2\,000  &  35.1  \\
1000-3D  &  3\,000  &  35.1  \\
1000-1D  &  10\,000 &  34.8  \\
\hline
Background & no data &  60.8 \\
\hline
\end{tabular}
\end{table}

\subsection{Comparison with the ``Eulerian'' method}\label{sec:43}
\begin{figure}
\begin{center}
\includegraphics[width=\textwidth]{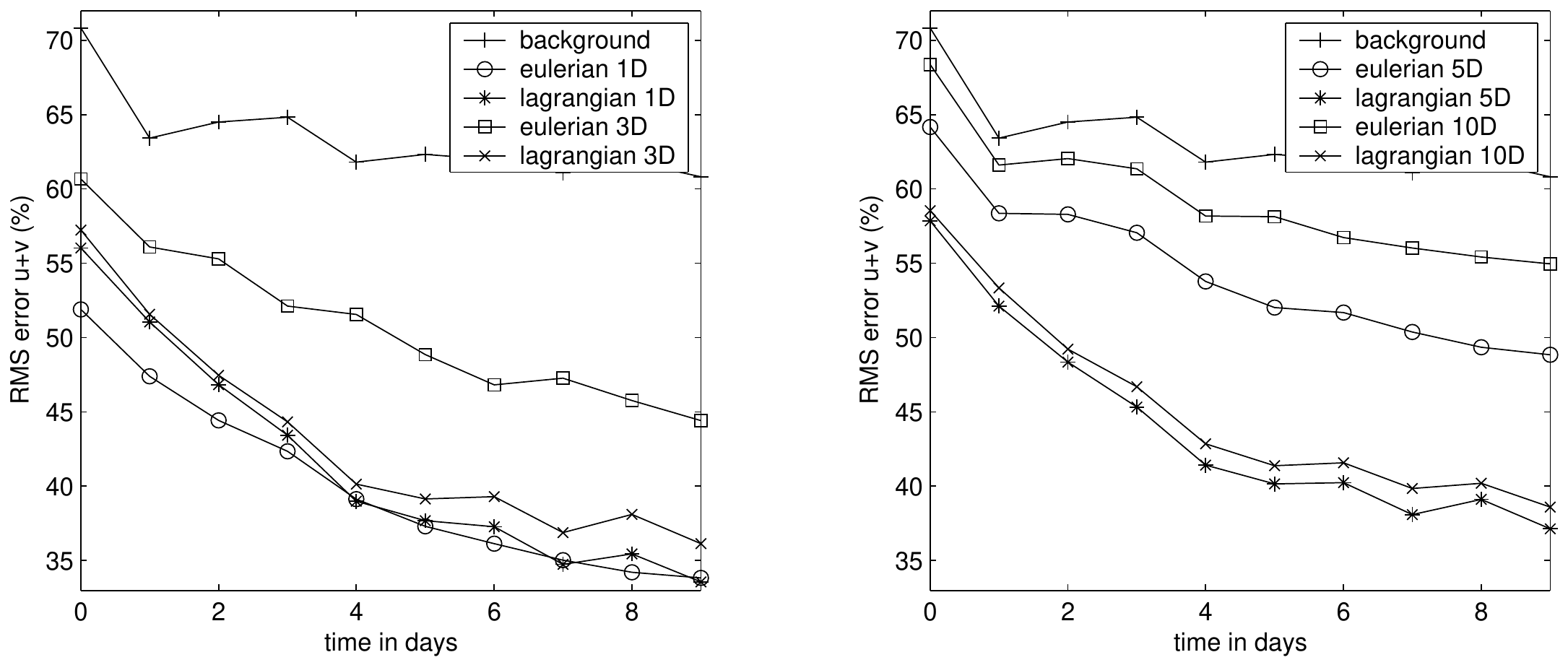}
\end{center}
\caption{\label{fig:10}Comparison Eulerian/Lagrangian: time-evolution of the RMS error for u+v corresponding to the assimilation of 3\,000 floats with different TSP. On the left, errors for the Lagrangian and Eulerian methods with small TSP. On the right, errors for the Lagrangian and Eulerian methods with large TSP. For reference, the background error is also displayed on each plot.}
\end{figure}

\begin{figure}
\begin{center}
\includegraphics[width=\textwidth]{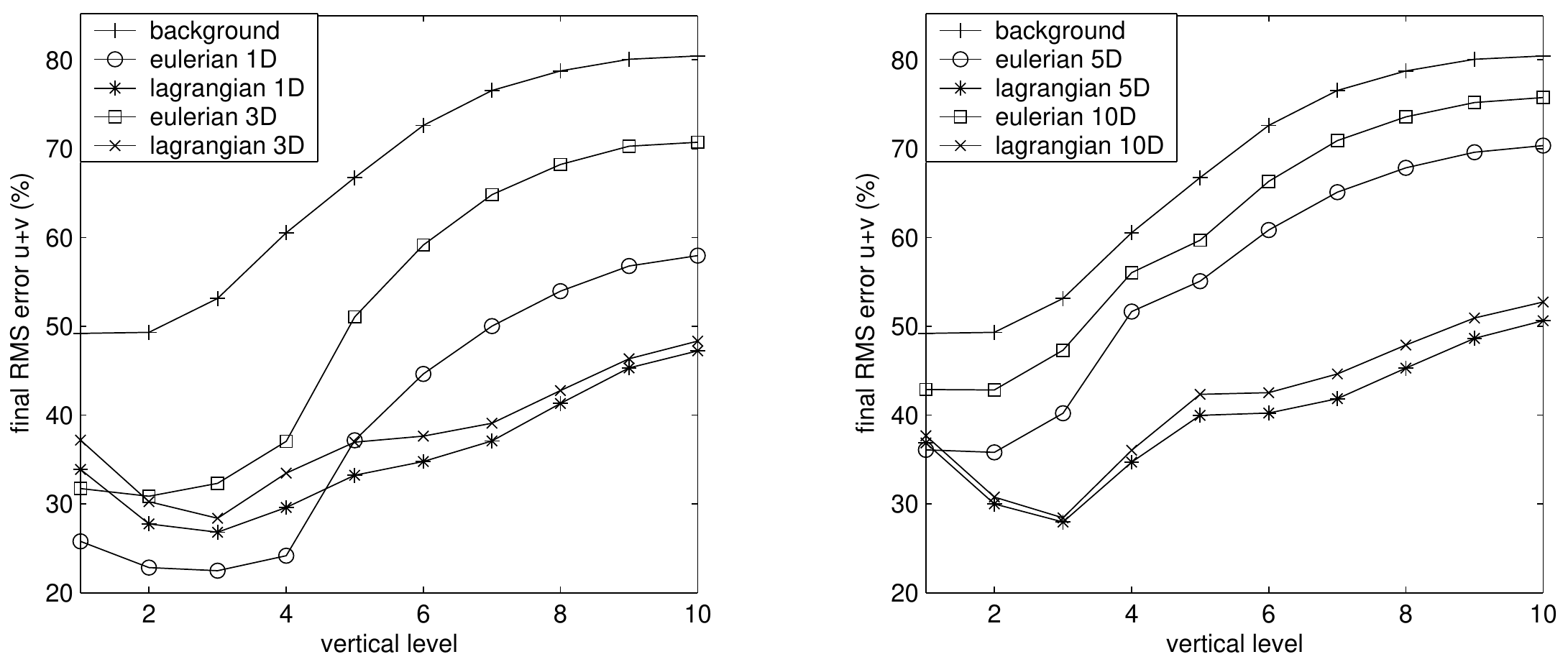}
\end{center}
\caption{\label{fig:11}Comparison Eulerian/Lagrangian: final RMS error for u+v as a function of the vertical level, corresponding to the assimilation of 3\,000 floats with different TSP. On the left, errors for the Lagrangian and Eulerian methods with small TSP. On the right, errors for the Lagrangian and Eulerian methods with large TSP. For reference, the background error is also displayed on each plot.}
\end{figure}

A classical method in oceanography is to assimilate the velocity observations deduced from the Lagrangian data according to the following finite differences formula:
\be 
\ba{lll}
\label{eq:10}
\displaystyle \frac{\xi_1(t_{k+1})-\xi_1(t_k)}{t_{k+1}-t_k} &\approx& u(\xi_1(t_k),\xi_2(t_k),z_0,t_k)\\
\displaystyle \frac{\xi_2(t_{k+1})-\xi_2(t_k)}{t_{k+1}-t_k} &\approx& v(\xi_1(t_k),\xi_2(t_k),z_0,t_k)
\ea 
\ee
Then the velocity data are treated as Eulerian data (measured at non-fixed points). We implement this method in the 4D-Var framework. The observation operator is much easier to write (and to differentiate and transpose) because it is an interpolation at the points of the true (fixed) floats trajectories. We compare the results for this method said ``Eulerian'' and for our ``Lagrangian'' one. Experiments have the same characteristics (3\,000 floats and varying TSP), their names are LAG-xxx or EUL-xxx with xxx the TSP, where xxx is the TSP. \\
Figures \ref{fig:10} and \ref{fig:11} represent the RMS error as a function of time and vertical level. We can see that the ``Eulerian'' approach is slightly better than the ``Lagrangian'' one when the TSP is small (one day), moreover its computation time is 10\% lower. Indeed the TSP is small enough so that the formula (\ref{eq:10}) is a very good approximation: the displacement vector between two successive positions is quasi tangent to the trajectory (ie quasi collinear with the velocity vector). However we can see that the error for the ``Lagrangian'' method is more homogeneous as a function of the vertical level. For larger TSP (3 days or more) the ``Lagrangian'' method is obviously better than the ``Eulerian'' one: the approximation formula (\ref{eq:10}) is not valid any more. We can see that our method is able to extract information from the positions data even if the TSP is large.

\subsection{Assimilation of noisy observations}\label{sec:44}
\begin{figure}
\begin{center}
\includegraphics[width=\textwidth]{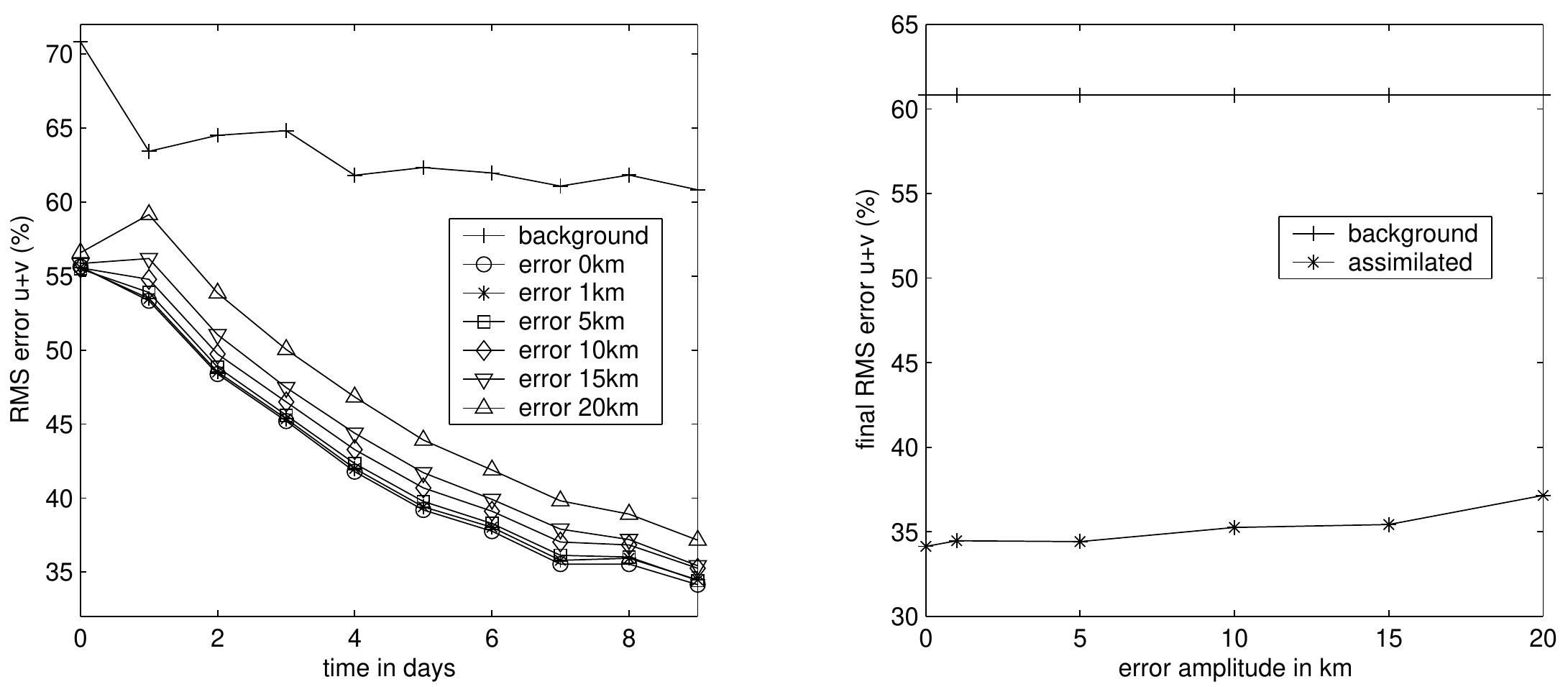}
\end{center}
\caption{\label{fig:12}Impact of observation errors: u+v RMS errors corresponding to assimilation of noisy observations. On the left: error as a function of time for experiments with noise amplitude from 1km to 20km. For reference, results for experiments without noise (0km) and without assimilation (background) are also displayed. On the right: final error as a function of the amplitude of the noise; the error without assimilation is also displayed.}
\end{figure}

\begin{figure}
\begin{center}
\includegraphics[width=\textwidth]{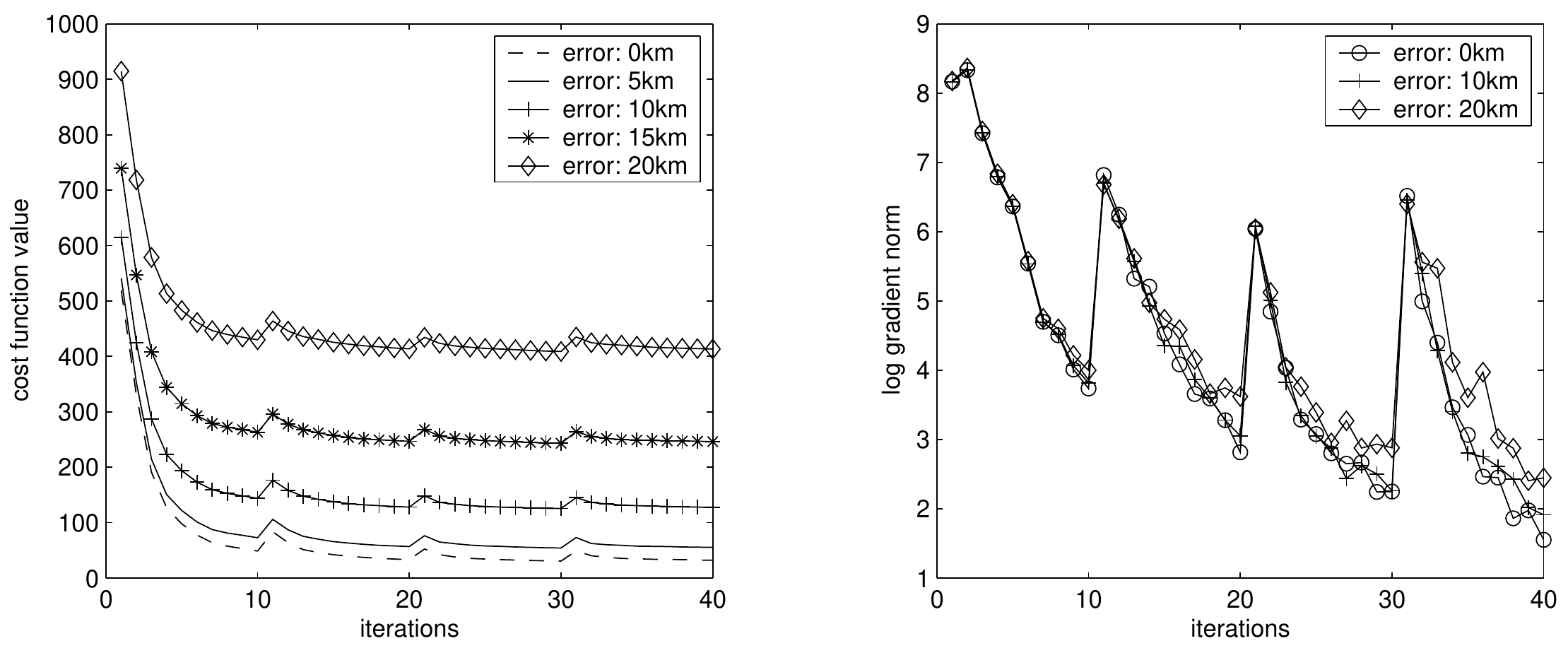}
\end{center}
\caption{\label{fig:13}Evolution of the cost function and its gradient's norm during the assimilation of noisy data. On the left, evolution of the cost function for experiments with noise amplitude from 1km to 20km. For reference, the cost for experiment without noise (0km) is also plotted. On the right, evolution of the gradient norm, shown on a logarithmic scale, for experiments with and without noise in data.}
\end{figure}

In order to deal with real data issues, a necessary first-step is to study the impact of observation errors in the twin experiments framework. To do that, we simulate as previously perfect data from the ``true state'' with 1\,000 floats drifting at level 4, their positions being sampled once a day. Then we add a random Gaussian noise to the computed positions. In the sequel, the word ``error'' represents the amplitude of the noise. As we said before, real errors are about 2 to 6 kilometers. However our system is idealized so we study the impact of errors up to 20 kilometers. The total displacement of one float between initial and final positions is around 25 kilometers (in steady regions) to 90 km (in the mid-latitude jet region), so that a 10 to 20-km noise is significant for most of the floats. \\
Figure \ref{fig:12} represents RMS errors as a function of time (on the left) for experiments with 0 to 20km errors and for the background (without assimilation). On the right we plot RMS errors as a function of observation error amplitude. The RMS error is very stable with increasing noise amplitude: our method is able to extract information even when the error amplitude is not negligible with respect to floats displacement.\\
Figure \ref{fig:13} shows the evolution of the cost function value and its gradient during the assimilation process. The abscissa represents the total number of iterations. Here we perform four outer loops: in each outer loop we perform ten inner minimization loops, in which the cost function is minimized, as we can see on the left. We can see that the gradient norm decreases and the cost function converges even in the presence of noise in data.

\section{Conclusion} \label{sec:5}
This paper shows that the problem of assimilating Lagrangian data can be solved by a variational adjoint method into a realistic primitive equations ocean model. We have implemented a Lagrangian method which takes into account the four dimensional (space and time) nature of the observations: RMS errors with assimilation are twice lower than without and the main patterns of the fluid flow are well identified at each vertical level, although the floats drift at a single determined level.\\
We have tested the sensitivity of our method to the characteristics of the dataset. It is very sensitive to the vertical drift level, and the best results are obtained for intermediates ones, especially level 4 (around 1000 meters depth). It is also very sensitive to the number of floats, but the more is not the better, it seems useless to launch more than 1\,000 floats in our configuration. It is very robust with respect to the increase of the time-sampling period, up to ten days.\\
We have compared our Lagrangian method to the Eulerian one, which consists in interpreting Lagrangian data as velocity information. When the time-sampling period of the observations is one day or less, the Eulerian method performs slightly better, but the transfer of information to lower levels is better achieved by the Lagrangian one. When this period is larger than two or three days, the Lagrangian method performs much better than the Eulerian one.\\
We also studied the impact of errors on observation: the reconstruction of the velocities is well achieved even with a large noise in data.\\
Also the performances of this method have been assessed in the framework of the twin experiments approach. The next step would be to use real data and to deal with problems such as trajectories modelling and model error.

\section*{Acknowledgments}
The author thanks Anthony Weaver and LODYC for the source code of OPAVAR 8.1.\\ Numerical computations were performed on the NEC SX5 vector computer at IDRIS.\\ This work is supported by French project Mercator.



\end{document}